\documentclass[12pt]{article}
\usepackage{amssymb}
\usepackage{amsthm}
\usepackage{amsxtra}
\usepackage{amsfonts}
\usepackage{xr}

\begin{document}

\def\sect{\section}

\newtheorem{thm}{Theorem}[section]
\newtheorem{cor}[thm]{Corollary}
\newtheorem{lem}[thm]{Lemma}
\newtheorem{prop}[thm]{Proposition}
\newtheorem{propconstr}[thm]{Proposition-Construction}

\theoremstyle{definition}
\newtheorem{para}[thm]{}
\newtheorem{ax}[thm]{Axiom}
\newtheorem{conj}[thm]{Conjecture}
\newtheorem{defn}[thm]{Definition}
\newtheorem{notation}[thm]{Notation}
\newtheorem{rems}[thm]{Remarks}
\newtheorem{rem}[thm]{Remark}
\newtheorem{question}[thm]{Question}
\newtheorem{example}[thm]{Example}
\newtheorem{problem}[thm]{Problem}
\newtheorem{excercise}[thm]{Exercise}
\newtheorem{ex}[thm]{Exercise}

\def\Bbb{\mathbb}
\def\cal{\mathcal}
\def\mL{{\mathcal L}}
\def\mC{{\mathcal C}}

\overfullrule=0pt

\def\si{\sigma}
\def\prf{\smallskip\noindent{\it        Proof}. }
\def\call{{\cal L}}
\def\nat{{\Bbb  N}}
\def\la{\langle}
\def\ra{\rangle}
\def\inv{^{-1}}
\def\ld{{\rm    ld}}
\newcommand{\ild}{{\rm ild}}
\def\trdeg{{\rm tr.deg}}
\def\dim{{\rm   dim}}
\def\th{{\rm    Th}}
\def\rest{{\lower       .25     em      \hbox{$\vert$}}}
\def\ch{{\rm    char}}
\def\zee{{\Bbb  Z}}
\def\conc{^\frown}
\def\acl{{\rm acl}}
\def\cls{cl_\si}
\def\cals{{\cal S}}
\def\mult{{\rm  Mult}}
\def\calv{{\cal V}}
\def\aut{{\rm   Aut}}
\def\ffi{{\Bbb  F}}
\def\ffiti{\tilde{\Bbb          F}}
\def\degs{deg_\si}
\def\calx{{\cal X}}
\def\gal{{\cal G}al}
\def\cl{{\rm cl}}
\def\loc{{\rm locus}}
\def\calg{{\cal G}}
\def\calq{{\cal Q}}
\def\calr{{\cal R}}
\def\caly{{\cal Y}}
\def\aff{{\Bbb A}}
\def\cali{{\cal I}}
\def\calu{{\cal U}}
\def\epsilon{\varepsilon} 
\def\Uu{{\cal U}}
\def\rat{{\Bbb Q}}
\def\ga{{\Bbb G}_a}
\def\gm{{\Bbb G}_m}
\def\cee{{\Bbb C}}
\def\ree{{\Bbb R}}
\def\frob{{\rm Frob}}
\def\Frob{{\rm Frob}}
\def\fix{{\rm Fix}}
\def\Uu{{\cal U}}
\def\proj{{\Bbb P}}
\def\sym{{\rm Sym}}
 
\def\dcl{{\rm dcl}}
\def\calm{{\mathcal M}}
\newcommand{\SU}{{\rm SU}}
\newcommand{\Cb}{{\rm Cb}}
\newcommand{\acb}{\overline{\rm Cb}}

\font\helpp=cmsy5
\def\semdp
{\hbox{$\times\kern-.23em\lower-.1em\hbox{\helpp\char'152}$}\,}

\def\dnfo{\,\raise.2em\hbox{$\,\mathrel|\kern-.9em\lower.35em\hbox{$\smile$}
$}}
\def\dnf#1{\lower1em\hbox{$\buildrel\dnfo\over{\scriptstyle #1}$}}
\def\dfo{\;\raise.2em\hbox{$\mathrel|\kern-.9em\lower.35em\hbox{$\smile$}
\kern-.7em\hbox{\char'57}$}\;}
\def\df#1{\lower1em\hbox{$\buildrel\dfo\over{\scriptstyle #1}$}}        
\def\stab{{\rm Stab}}
\def\qfcb{\hbox{qf-Cb}}
\def\perf{^{\rm perf}}
\def\sipm{\si^{\pm 1}}

\def\vlabel{\label}

\title{Remarks around the non-existence of difference-closure}

\author{Zo\'e Chatzidakis\thanks{partially supported  by
    ANR--13-BS01-0006 (ValCoMo). }\\ DMA (UMR 8553), Ecole Normale Sup\'erieure\\
CNRS, PSL Research University}

\maketitle

\begin{abstract} This paper shows that in general, difference fields do not have a difference closure. However, we introduce a stronger notion of closure ($\kappa$-closure), and show that every algebraically closed difference field $K$ of characteristic 0, with fixed field satisfying a certain natural condition, has a $\kappa$-closure, and this closure is unique up to isomorphism over K. 
  \end{abstract}
\section*{Introduction}
In this paper,  a difference field is a field $K$ with a
distinguished automorphism $\si$. A difference  field $L$ is
{\em difference-closed} if every finite system of difference equations with
coefficients in $L$ which has a solution in a difference field extending
$L$, has already a solution in $L$. \\
If $K$ is a difference field, then a {\em difference-closure of $K$} is
a difference-closed field containing $K$, and which $K$-embeds into
every difference-closed field containing $K$. \\[0.05in]
The algebra of difference fields
was developped by Ritt, in analogy with the algebra of differential
fields. It is well-known that any differential field of characteristic
$0$ has a differential closure, and that this differential closure is
unique up to isomorphism over the field. In 2016, Michael Singer asked whether  this result 
generalises to the context of difference fields. One of the main results
of this paper is that 
it does not, even after imposing some natural conditions on the 
difference field $K$. We will show by two examples (\ref{ex1} and
\ref{ex2}) that even the
existence of a difference closure can fail.  \\[0.05in]
There are several natural
strengthenings 
of the notions of difference-closed and difference-closure (originating from model-theory but
which have a natural algebraic translation), and we will show that these notions do satisfy
existence and uniqueness of  closure, provided we work over an
algebraically closed difference
field of {\bf characteristic $0$} whose fixed subfield is large
enough. \\[0.05in]
The theory of difference-closed difference fields has been extensively
studied, and is commonly denoted by ACFA. The proof of our result uses
in an essential way the characteristic $0$ hypothesis, as it allows us 
to use techniques
of stability theory. They  provide  examples of  structures
which are stable over a predicate, see \cite{PS} and \cite{SU} for
definitions.  The  main result of the paper is\\[0.05in]
{\bf Theorem \ref{thm1}}. {\em Let $\kappa$ be an
uncountable cardinal or $\aleph_\epsilon$, and let $K$ be an algebraically closed difference field of
characteristic $0$ such that $F:=\fix(\si)(K)$ is pseudo-finite and is
$\kappa$-saturated. Then there is a $\kappa$-prime model of ACFA 
over $K$. Furthermore, it
is unique up to isomorphism over $K$.  }\\[0.05in]
Here is  an algebraic translation of this result for $\kappa\geq
\aleph_1$: 
Call a difference field $\calu$ {\em $\kappa$-closed} if every system of $<\kappa$
difference equations over $\calu$ which has a solution in some difference
field extending $\calu$, has a solution in $\calu$. 
The field $\calu$ is a
{\em $\kappa$-closure} of the difference field $K$ if it is $\kappa$-closed,
contains $K$, 
and $K$-embeds into every $\kappa$-closed difference field containing
$K$. Then Theorem \ref{thm1} states, for $\kappa\geq \aleph_1$:\\[0.05in]
{\em Let $K$ be an algebraically closed difference field of
  characteristic $0$, whose fixed field  $F$ is pseudo-finite and such that
  every system of $<\kappa$ polynomial equations over $F$ which has a solution in
  a regular extension of $F$, already has a solution in $F$. 
  Then $K$ has a $\kappa$-closure, and it is unique
  up to $K$-isomorphism.}\\[0.1in]
It is unlikely that this result can be generalised to the characteristic
$p$ context, and in fact, I conjecture that unless the difference field
$K$ of characteristic $p>0$ is of cardinality $<\kappa$ or is already $\kappa$-closed, then it
does not have a $\kappa$-closure. \\[0.1in]
The paper is organised as follows. In section 1 we discuss the problem
and 
reformulate it in model-theoretic terms, and describe the two
examples. In section 2, we state the preliminary results we will need
from difference algebra and model theory. Section 3 contains the proof
of Theorem~\ref{thm1}.

\section{Discussion of the problems and the examples}

\para{\bf Notation and conventions}. \vlabel{notation} All difference fields will be inversive, i.e., the endomorphism $\si$
will be onto. Let $K$ be a difference field, contained in some large
difference field $\calu$. If $a$ is a tuple in $\calu$, we denote by
$K(a)_\si$ the difference field generated by $a$ over $K$, i.e., the
subfield $K(\si^i(a))_{i\in\zee}$ of $\calu$. The algebraic and
separable closure of a field $L$ are denoted by $L^{alg}$ and $L^s$
respectively, and $G(K)$ denotes the absolute Galois group of $K$, i.e.,
$\gal(K^s/K)$. If  $A\subset \calu$, then $\acl(A)$ denotes the
smallest algebraically closed difference field containing $A$; it
coincides with the model-theoretic algebraic closure of $A$ for the
theory ACFA (1.7 in \cite{CH}). We denote
by $\call$ the language $\{+,-,\cdot,0,1,\si\}$.

\para {\bf Translation into model-theoretic terms}. Let $K$ be a
difference field. 
Recall that any complete theory extending the theory ACFA of difference-closed difference fields
is supersimple, unstable, of SU-rank $\omega$, and does not eliminate
quantifiers, but it eliminates imaginaries. It is extensively studied in
\cite{CH}. The reason ACFA does not eliminate quantifiers is that given
an automorphism $\si$ of a field $K$, there may be several
non-isomorphic ways of extending $\si$ to $K^{alg}$. So, the first
obvious obstacle to the existence of a difference closure is that, and
a natural condition to impose is to assume that $K$ is algebraically
closed. There is another natural condition one needs to impose: if $L$ is difference closed, then its fixed field $\fix(\si)(L)=\{a\in
L\mid \si(a)=a\}$ is pseudo-finite. Moreover, every pseudo-finite field
can occur as the fixed field of some difference closed field (\cite{A}). Thus if $L$ is the
difference-closure of a difference field $K$, then $\fix(L)$ must be
prime over $\fix(K)$ (for its theory in the language of rings). 
Duret showed in
\cite{Du} that any completion of the theory of pseudo-finite fields has
the independence property. From his proof one extracts easily the fact
that non-algebraic types are non-isolated, and this forces us to require
in case $K$ is countable that $\fix(\si)(K)$ be
pseudo-finite in order to hope to have a difference closure. The case when $K$ is
uncountable is a little more complicated, the question is adressed and
solved in \cite{C-Psf}. \\[0.05in]
It is therefore 
reasonable 
to make the following two assumptions: 
\begin{quote}{\em $K$ is algebraically closed, and $K\cap \fix(\si)$ is
  pseudo-finite}.
\end{quote}
But even this is not enough. To show this does not suffice, what we
need to do is the following:
\begin{quote}{\em Exhibit a difference field $K$ satisfying the above
 two   conditions, and a finite system of difference equations which does
    not have a solution in $K$, and such that any finite strengthening
    of this system has several completions.} 
\end{quote}
This looks easy, since even our stable types are only supersimple, not
$\omega$-stable. However, the first obvious examples \dots\ do not
satisfy the first condition. Here is a more involved example, taken from
\cite{CH} (example 6.7): 

\begin{example}{\vlabel{ex1}} An example in characteristic $0$.\\
Let $k$ be a 
  pseudo-finite field of characteristic $0$ containing $\rat^{alg}$, and
  consider $K=(k^{alg},\si)$, where $\si$ is a (topological) 
  generator of $\gal(K/k)$. We conside the elliptic curve $J$, with
  $j$-invariant $a\notin  K$, and which is defined by 
$$y^2+xy=x^3-\frac{36}{a - 1728}x -\frac{1}{a-1728}.$$
We let $A'$ be a cyclic subgroup of $J$ of order $p^2$, $A=[p]A'$ and $a_1$ the
$j$-invariant of the elliptic curve $J/A$, $a_2$ the $j$-invariant of
the elliptic curve $J/A'$. Then the map $\rat^{alg}(a,a_1)\to \rat^{alg}(a_1,a_2)$
which is the identity on $\rat^{alg}$ and sends $(a,a_1)$ to $(a_1,a_2)$  extends
to a field automorphism of $\rat(a)^{alg}$, which in turns extends to an
automorphism of $K(a)^{alg}$ which agrees with $\si$ on $K$. Let
$\Phi(x,x_1,x_2)$ be the finite system of polynomial equations  which describe
the algebraic locus of $(a,a_1,a_2)$ over $K$, and consider the difference equation
$$\Phi(x,\si(x),\si^2(x)).$$
The formula $\Phi$ can be expressed, in the notation of Chapter 5 \S3 of
\cite{La-Ell}  (see in particular Theorem
5), 
 by $\Phi_{p^2}(\si^2(x),x)=
\Phi_{p}(\si^2(x),\si(x))=\Phi_p(\si(x),x)=0$. 
Let $b$ be any solution of $\Phi$ which satisfies $\si(b)\neq b$. Note
that necessarily, the kernel of the map $J_b\to J_{\si^n(b)}$ for $n>0$ is
cyclic of order $p^n$. As $\si(b)\neq b$, we know that $b$ is
transcendental. Hence the curve $J_b$ is not of CM-type, its
endomorphism group is isomorphic to $\zee$, and therefore $J_b$ is not
isomorphic to any of its quotients by  finite cyclic
subgroups (see e.g. section C.11 in \cite{Sil}). Therefore,   the elements $b,\si(b),\si^2(b), \ldots$
are all distinct, and $b\notin K$. Furthermore the isomorphism type of
$K(b)_\si$ over $K$ is determined by $\Phi(b,\si(b),\si^2(b))$. \\
So any difference-closed field containing $K$ must contain a solution of
$\Phi(x,\si(x),\si^2(x))\land \si(x)\neq x$. However, Example 6.7 of
\cite{CH} shows that if $b$
is as above, and $L$ is any finite extension of $K(b)_\si$, then there
are $2^{\aleph_0}$ non-isomorphic ways of extending $\si$ to
$L^{alg}$. Thus $K$ does not have a difference closure.

\end{example}
\noindent
One can build other examples along the same lines, using moduli spaces
of abelian varieties. 

\begin{example}\vlabel{ex2} An example in characteristic $p>0$.\\ 
Let $K=k(A)_\si^{alg}$, where $k$ is  a
pseudo-finite field fixed by $\si$, $\si$ restricts to a generator of
$\gal(k^{alg}/k)$, and $A$ is the set of solutions of the equation
$\si(x)^p-\si(x)+x^p=0$ (in some countable difference-closed
overfield). So, any difference-closed field containing $K$ will contain
the set $B$ of solutions of the equation $\si(x)-x^p+x=0$, an
infinite-dimensional $\ffi_p$-vector space. However, as was shown in
example 6.5 of \cite{CH}, there are $2^{|A|}$ ways of extending
$\si$ from $Kk(B)_\si^{alg}$ to $K(B)_\si^{alg}$: there is a definable
non-degenerate bilinear map $q:A\times B\to \ffi_p$, which can be chosen
totally arbitrarily.\\[0.1in]
In fact this example is part of a large family of examples: let $f$ and
$g$ be additive polynomials with coefficients in a difference field $K$,
and assume that the subgroup $A$ of $\ga$ defined by $f(x)=g(\si(x))$ is
locally modular. Then there is a definable subgroup $B$ of $\ga$, and a
definable non-degenerate bilinear map $A\times B\to \ffi_p$. As above,
there is no prime model over $K(B)_\si$. 
\end{example}

\noindent 
While we provided examples of difference fields not having a difference
closure, we did not provide a procedure which, given a difference field
which is not difference closed, exhibits a non-isolated type which needs
to be realised.  So, the following remains open:

\begin{question} {\em Are there any  difference fields
  which are not difference closed but admit a difference closure?} 
\end{question}


\section{Preliminaries}

\bigskip\noindent
{\bf \large Basic difference algebra}. 

\para Let $K\subset\calu$ be difference fields. If
$X=(X_1,\ldots,X_n)$, the ring $K[X]_\si=K[\si^i(X_j)]_{1\leq j\leq n,
  i\in\nat}$ is called the {\em $n$-fold difference polynomial ring}. A {\em
  difference equation} is an equation of the form $f(X)=0$ for some
$f(X)\in K[X]_\si$.\\[0.05in]
If $a$ is a finite tuple in $\calu$, and $L$ is a difference subfield of
$K(a)_\si$ containing $K$, then $L=K(b)_\si$ for some finite tuple $b$
(5.23.18 in \cite{Co}). \\[0.05in]
An element $a\in \calu$ is {\em transformally
algebraic} over $K$ if it satisfies some non-trivial difference equation
with parameters in $K$. Otherwise, it is {\em transformally
  transcendental} over $K$. A tuple $a$ is {\em transformally algebraic over $K$}
if all its elements are. A (maybe infinite) tuple of elements of $\calu$ is
transformally independent over $K$ if it does not satisfy any
non-trivial difference equation with coefficients in $K$. A
{\em transformal transcendence basis of $\calu$ over $K$} is a subset
$B$ of
  $\calu$ which is transformally independent over $K$ and maximal such;
  every element of $K$ will then be transformally algebraic over
  $K(B)_\si$. We denote by $\Delta(K)$, the {\em transformal transcendence
  degree of $K$}, i.e., the cardinality of a transformal transcendence basis
of $K$, and if $L$ is difference field containing $K$, by $\Delta(L/K)$
the cardinality of a transformal transcendence basis of $L$ over $K$. 


\para{\bf The fixed field}. \vlabel{fix1} The fixed field of $\calu$ is the field
$\fix(\si)(\calu):=\{a\in \calu\mid \si(a)=a\}$. Then $\fix(\si)(\calu)$
and $K$ are linearly disjoint over their intersection. (Choose $n$
minimal such that there are $c_1, c_2,\ldots,c_n\in\fix(\si)$ and
$d_1=1,d_2,\ldots,d_n\in K$ such that $\sum_ic_id_i=0$; applying $\si$ we get
$\sum c_i\si(d_i)=0$, and by minimality of $n$, that $\si(d_i)=d_i$ for
all $i$.)  This implies in particular that if $E$ is a difference
subfield of $K$, then $E\fix(\si)(\calu)$ and $K$ are linearly disjoint
over their intersection $E(\fix(\si)(\calu)\cap K)$. In positive
characteristic, similar results hold for the other fixed fields
$\fix(\si^n\frob^m)$.  \\[0.15in]

%
\noindent 
{\large \bf Basic model-theoretic facts}.

\para For references see \cite{CH}. The theory ACFA is supersimple, of SU-rank
$\omega$. It eliminates imaginaries, but does not eliminate quantifiers. The completions of ACFA are given by describing the
isomorphism type of the automorphism $\si$ of the algebraic closure of
the prime field (1.4 in \cite{CH}). \\
Let $\calu$ be a model of ACFA, and
$K$ a difference subfield of $\calu$. 

\para {\bf Types, algebraic closure, independence}. \label{types} If $a$ is a tuple of elements of $\calu$, then $tp(a/K)$ is determined by
the isomorphism type of the difference field $\acl(Ka)=K(a)_\si^{alg}$ over $K$:
$a$ and $b$ have the same type over $K$ if and only if there is a
$K$-isomorphism of difference fields $K(a)_\si\to K(b)_\si$ which send
$a$ to $b$ and extends
to the algebraic closure of $K(a)_\si$. (1.5 in
\cite{CH}). The SU-rank of $a$ over $K$, denoted by $\SU(a/K)$, is bounded
by $\trdeg(K(a)_\si/K)$, and is finite if and only if
$\trdeg(K(a)_\si/K)$ is finite (if and only if $a$ is transformally
algebraic over $K$). \\
Let $A,B,C$  be subsets of $\calu$. Then $A$ is {\em independent from $B$
over $C$}, denoted $A\dnfo_C B$, if and only if the fields $\acl(AC)$ and $\acl(BC)$ are
free over $\acl(C)$. Equivalently, if whenever $a$ is a tuple of
elements in $A$, then the prime $\si$-ideal $I_\si(a/\acl(BC)):=
\{f(X)  \in 
\acl(BC)[X]_\si\mid f(a)=0\}$ is generated (as a $\si$-ideal) by its intersection with
$\acl(C)[X]_\si$. Then independence coincides with non-forking, and we
will also say, in that case, that {\em $tp(A/B)$ does not fork over $C$}.

\para {\bf Reducts}. Let $n>0$ be an integer. We denote by $\call[n]$
the language $\{+,-,\cdot,0,1,\si^n\}$, and by $\calu[n ]$ the reduct
$(\calu, \si^n)$ to the language $\call[n]$. By Corollary 1.12 of
\cite{CH}, $\calu[n]\models {\rm ACFA}$. If $a$ is a tuple in $\calu$,
then $tp(a/K)[n]$ denotes the type of $a$ in the reduct $\calu[n]$, and
$qftp(a/K)[n]$ the quantifier-free type of  $a$ in the reduct
$\calu[n]$. 

\para{\bf Notions of canonical bases}. \vlabel{cb} If $a$ is a tuple in $\calu$,
then $\Cb(a/K)$ denotes the smallest difference field over which
$I_\si(a/K)$ is defined.  Then $tp(a/K)$ does not
fork over $\Cb(a/K)$. Also, $\Cb(a/K)$ is contained in the algebraic
closure over $K$ of $n$ independent realisations of $tp(a/K)$, for some $n$; if
$K(a)_\si$ is a regular extension of $K$, then $\Cb(a/K)$ is
contained in the  difference field generated over
$K$ by finitely many independent realisations of $tp(a/K)$ (see the proof of Lemma 2.13(4) in
\cite{CH}). 
$\acb(a/K)$ will denote $\Cb(a/K)^{alg}$. Note that a type does not fork over the
algebraic closure of some finite set. 
 
\para{\bf The generic type}. \vlabel{generic} The {\em generic} $1$-type is the type of a
transformally transcendental element. It is axiomatised by its
quantifier-free part, is definable and stationary\footnote{A type $p$
  over a set $A$ is stationary if whenever $B\supset A$, then $p$ has a
  unique non-forking extension to $B$.}. Similarly, if $V$ is
a variety defined over the algebraically closed difference field $K$,
then the generic type of $V$ (which is characterised by having a
realisation $a$ with $\Delta(K(a)_\si/K)=\dim(V)$) is axiomatised by its quantifier-free part,
is definable and stationary (see 2.11 in \cite{CH}). 

\para{\bf Orthogonality of types}.  Let $p$ and $q$ be (partial) types  over $A$ and
$B$ respectively. \\
If $A=B$, we say that $p$ and $q$ are {\em almost
  orthogonal} (or {\em weakly orthogonal}), denoted by $p\perp^a q$, if
whenever $a$ realises $p$ and $b$ realises $q$, then $a\dnfo_A b$. \\
We say that $p$ and $q$ are {\em orthogonal}, denoted by $p\perp q$, if
whenever $C$ contains $A\cup B$, and $a$ realises $p$, $b$ realises $q$,
and $a\dnfo_A C$, $b\dnfo_B C$, then $a\dnfo_C b$. 

\para{\bf The dichotomy in characteristic $0$}. \vlabel{dich} Recall that a partial type $\pi$ over a set
$A$ is called {\em one-based}\footnote{In \cite{CH}, they are called
  {\em modular}.} if whenever $a_1,\ldots,a_n$ realise $\pi$ and
$B\supset A$, then $(a_1\ldots a_n)\dnfo_C B$, where
$C=\acl(Aa_1,\ldots,a_n)\cap \acl(B)$\footnote{Here we are using the
  fact that any completion of ACFA eliminates imaginaries}. \\
Types of finite $\SU$-rank are
analysable in 
terms of types of $\SU$-rank $1$. The main result of \cite{CH} says that
in characteristic $0$, a type $q$ of SU-rank $1$ is either
one-based, or non-orthogonal to the fixed field. Moreover, if  $q$ is one-based, then it is stable
stably embedded, and definable. See Theorem 4.10 in \cite{CH}.

\para{\bf Stable embeddability of the fixed field}. \vlabel{fix2} Recall that a subset
$S$ of $\calu^n$, which is definable or $\infty$-definable, is {\em stably
embedded} if whenever $D\subset \calu^{nm}$ is definable with parameters
from $\calu$, then $D\cap
S^m$ is definable with parameters from $S$. An important result of
\cite{CH} (Proposition 1.11) says that the fixed field $F:=\fix(\si)$ of $\calu$ is
stably embedded: if $D\subset F^n$ is definable in the
difference field $\calu$ (with parameters from $\calu$), then it is definable in the pure field
language in $F$ (with parameters from
$F$).  In fact, one has more: let $C=\acl(C)\subset
\calu$, and $b$ a tuple in $F$; then $tp_F(b/C\cap F)\vdash
tp_{\calu}(b/C)$: indeed, all finite $\si$-stable extensions of
$CF$ are contained in $CF^{alg}$ (see Lemma 4.2 in
\cite{CHS}), and therefore any $(C\cap F)$-automorphism of the field $F$
extends to a $C$-automorphism of the difference field $\acl(CF)$,
since it obviously extends to a $C$-automorphism of $CF$, and the
automorphism $\si$ of $CF^{alg}$ extends uniquely to $\acl(CF)$
(Babbitt's Theorem. See e.g. Lemma 2.8 in \cite{CH}). \\
For more properties of stably embedded sets or types,
see the appendix of \cite{CH}. 

\para\vlabel{sse} {\bf More on stable stably embedded types}. For a
definition of a (partial) type being {\em stable stably embedded}, see Lemma 2
of the appendix of \cite{CH}. Here we will use the following
consequence: let $A=\acl(A)$ be algebraically closed, and suppose that
$tp(a/A)$ is stable stably embedded. Then $tp(a/A)$ is definable (over
$A$; see Lemma 1 in the Appendix of \cite{CH}). Also, if $B=\acb(a/A)$ and $tp(a/B)\perp^a tp(A/B)$, then
$tp(a/B)\vdash 
tp(a/A)$: this is because $tp(a/B)$ has a unique non-forking extension
to any superset of the algebraically closed set $B$.

\begin{defn} \vlabel{internality}({\rm \bf Internality to the fixed
    field}). Let $\pi$  be a partial type 
over $A\subset \calu$, and $F=\fix(\si)(\calu)$. 
\begin{enumerate} 
\item $\pi$ is qf-internal to $F$ if there is some finitely generated
  over $A$ difference field $C$
 such that whenever $a$ realises $p$, there is a tuple
  $b$ in $F$ such that $a\in C(b)$. I.e., $a\in CF$. 
\item $\pi$ is almost-internal to $F$ if there is some some finitely generated
  over $A$ difference field $C$ such that whenever $a$ realises $p$, there is a tuple
  $b$ in $F$ such that $a\in \acl(Cb)$. 
\end{enumerate}
\end{defn}
\begin{rems} \vlabel{internality-2}
  Clearly qf-internality implies almost-internality. Moreover, to show
  qf-internality or almost-internality of a (complete) type $p$, it is enough to do it for a
  particular realisation $a$ of the type $p$, i.e., to find $C$
  independent from $a$ over $A$ such that $a\in CF$ or $a\in
  \acl(CF)$. See Lemma 5.2  in \cite{CH}. \\
  Internality or almost internality (to $F$) of a type is in fact a property of
  its quantifier-free part.\\
Recall that
a difference field $E$ is linearly disjoint from $F$ over
$F\cap E$. It follows that 
in (1) above, the tuple $b$ can be taken so that $C(b)=C(a)_\si$: take a
generating  tuple $d$ of $F\cap C(a)_\si$ over $F\cap C$; as $F$ is linearly
disjoint from $C(a)_\si$ over $F\cap C(a)_\si$, we get that $C(a)_\si=C(d)$. 
\end{rems}
\begin{lem}\vlabel{internality-3} Let $A=\acl(A)$, and assume that
  $tp(a/A)$ is almost internal to $F$. Then there is $a'\in A(a)_\si$
  such that $tp(a'/A)$ is qf-internal to $F$, $\si(a')\in A(a')$, and $a\in \acl(Aa')$. 
\end{lem}

\prf  By assumption there is
some tuple $c$ independent from $a$ over $A$ and such that
$a\in \acl(AFc)$. Taking $b$ in $F$ such that
$A(c,a)_\si\cap F=(F\cap A)(b)$, we obtain that $F$ is linearly disjoint
from $A(c,a)_\si$ over $(F\cap A)(b)$, and therefore that $AF(c,b)_\si$
and $A(c,a)_\si$ are linearly disjoint over $A(c,b)_\si$, so that
$a\in\acl(Acb)$ (since $a\in \acl(AFcb)$).  As $c$ is independent  from $a$ over
$A=\acl(A)$, it follows that $A(c,a)_\si=A(c,a,b)_\si$ is a regular extension of
$A(a)_\si$, and therefore that $\Cb(b,c/A(a)_\si)$ is contained in the
difference field generated by finitely many realisations of
$tp(b,c/A(a)_\si)$ (see \ref{cb}). Again, as $c$ is independent from $a$ over $A$ and
$b$ is in $F$, it follows that if $a'$ is such that
$\Cb(c,b/A(a)_\si)=A(a')_\si$, then $tp(a'/A)$ is qf-internal to $F$. As
$b\in A(a',c)_\si$ and $c$ is independent from $a$ over $A$, it follows that
$a\in \acl(Aa')$ as desired. As $A(c,a')_\si =A(c,b)_\si$ and $b\in F$, it follows
that $A(c,a')_\si$ is finitely generated as a field extension of
$A(c)_\si$. But as $a'$ and $c$ are independent over $A$, the same holds
of the field extension $A(a')_\si/A$, i.e., for some $n$, $\si^n(a')\in
A(a',\si(a')\ldots,\si^{n-1}(a'))$. We then replace  $a'$ by 
$(a',\si(a')\ldots,\si^{n-1}(a'))$

\para{\bf The semi-minimal analysis}. \vlabel{semi-min} Let $a$ be a
tuple which is 
transformally algebraic over $K$. Thus $\SU(a/K)<\omega$. As Th$(\calu)$
is supersimple, there is a sequence $a_1,\ldots,a_n\in \acl(Ka)$, such
that $a\in\acl(Ka_1,\ldots,a_n)$, and for every $0<i\leq n$,
$tp(a_{i}/\acl(K,a_1,\ldots,a_{i-1}))$ is either one-based of rank $1$, or
almost-internal to a non-one-based type of rank $1$. This is a classical 
result in supersimple theories, for a proof in our case in
characteristic $0$, see Theorem~5.5 in \cite{CH}. Note that in
characteristic $0$ by the above dichotomy 
\ref{dich}, all non-one-based types of rank $1$ are non-orthogonal to
$\si(x)=x$,  and by Lemma \ref{internality-3}, almost-internality to $F$ may be replaced by
qf-internality to $F$. 
\begin{defn} \vlabel{epsilon1} Let $T$ be a completion of ACFA, $M$ a model of $T$. 
\begin{enumerate}
\item We say that $M$ is {\em $\aleph_\epsilon$-saturated} if whenever
    $A\subseteq M$ is finite, then every strong $1$-type over $A$ is
    realised in $M$. Equivalently, as our theory eliminates imaginaries,
    if every $1$-type over $\acl(A)$ is realised in $M$. 
\item Let $\kappa$ be an infinite cardinal or $\aleph_\epsilon$, and
  $A\subseteq M$. We say
  that $M$ is {\em $\kappa$-prime over $A$}, if $M$ is $\kappa$-saturated, and
  $A$-embeds elementarily into every $\kappa$-saturated model of
  Th$(M,a)_{a\in A}$. When $\kappa=\aleph_\epsilon$, one also says that
  $M$ is {\em
    a-prime}.

\item Let $\kappa$ be an infinite cardinal or $\aleph_\epsilon$. We say
  that  $A\subseteq M$ is {\em small} if  $A=\acl(A_0)$, where $A_0$ is finite if we are dealing
  with $\aleph_\epsilon$-saturation, and has cardinality $< \kappa$ if we
  are dealing with $\kappa$-saturation. 
We will also say that $A\subseteq M$ is {\em very small} if
$A=\acl(A_0)$, where $A_0$ is finite.
  \item Let $\kappa$ be an infinite cardinal or $\aleph_\epsilon$, and
  $A\subseteq M$. A type $p$ over $A$ is {\em $\kappa$-isolated} if it is
    implied by its restriction to some small subset of $\acl(A)$. 
\item We say
    that $M$ is {\em  $\kappa$-atomic} over $A\subseteq M$ if whenever $a$ is a
    (finite) 
    tuple in $M$, then $tp(a/A)$  is
    $\kappa$-isolated. Recall also that $M$ is atomic over $A$ if every
    finite tuple realises an isolated type over $A$. 
    \item We say that  $B=\acl(B)\subseteq M$ is {\em $\kappa$-constructed over $A\subseteq
      M$} if there is  a sequence
$(d_\alpha)_{\alpha<\mu}$ in $B\setminus A$ such that for every
$\alpha<\mu$, $tp(d_\alpha/\acl(Ad_\beta\mid \beta<\alpha))$ is
$\kappa$-isolated and $B=\acl(Ad_\alpha\mid \alpha<\mu)$. 
\end{enumerate}
\end{defn}

\begin{rems}\vlabel{atom} \begin{enumerate}
  \item If $\kappa$ is a regular cardinal, then $\kappa$-atomicity
  is transitive: if $A\subseteq B\subseteq C\subseteq M$, with $B$
  $\kappa$-atomic over $A$ and $C$ $\kappa$-atomic over $B$, then $C$ is
  $\kappa$-atomic over $A$. This is however not necessarily true when $\kappa$ is
  singular. However, this will hold if $B=\acl(Ab)$ for some finite
  tuple $b$ (since every finite tuple in $B$ realises an isolated type over
  $Ab$), or if $C$ is atomic over $B$. (There are stronger statements
  involving cardinals $\lambda<{\rm cf}(\kappa)$.)
  \item 
 If $M$ is a $\kappa$-saturated model of $T$ containing $A$ and $M$ is
 $\kappa$-constructed over $A$, then $M$ is $\kappa$-prime over
 $A$.
 \item The property of being $\kappa$-constructed is preserved under
   towers and union of chains indexed by ordinals.
   \end{enumerate}
  \end{rems}

\para {\large \bf Algebraic translation of the model-theoretic notions}\\[0.05in]
  Let us translate what the notions of saturation mean in our case. We
  will be dealing with either uncountable cardinals or
  $\aleph_\epsilon$. Recall that $tp(a/A)$ is entirely determined by the
  isomorphism type over the difference field generated by $A$ of  the
  difference field $\acl(Aa)$ (\ref{types}). So, for $\kappa$ an uncountable cardinal,
  the $\kappa$-saturation of a model $M$ of ACFA simply means: every system of $<\kappa$
  difference equations with coefficients in $M$, which has a solution in
  a difference field extending $M$, already has a solution in
  $M$. This is what was called $\kappa$-closed in the
  introduction. \\
The notion of $\kappa$-prime over a difference subfield corresponds to
being a $\kappa$-closure of that difference field.  \\
In the case of $\aleph_\epsilon$-saturation, the
  algebraic description is a little more complicated, and is better
  expressed in terms of embedding problems: 
    work inside a large model
  $\calu$, and consider a submodel $M$ of $\calu$. Then $M$ is
   $\aleph_\epsilon$-saturated if whenever $a$ is a finite tuple of
   elements of $M$ and $b$ an element of $\calu$, there is an
   $\acl(a)$-embedding of $\acl(a,b)$ inside $M$. \\
A similar description
   holds for $\kappa$-saturated, with the base set $a$ of cardinality
   $<\kappa$: a model $M$ of ACFA is $\kappa$ saturated if whenever
   $A\subset M$ is small and $b$ is a finite tuple in some difference field $\calu$
   containing $M$, then there is an $A$-embedding of $\acl(Ab)$ into
   $M$. Note that  $|A|$-many difference equations are necessary to
   describe the isomorphism type of $\acl(Ab)$ over $A$. 


\section{The results}

Results of Hrushovski (\cite{Hr-PAC}) show that if $F$ is a pseudo-finite
field and $C\subset F$, then ${\rm Th}(F,c)_{c\in C}$ eliminates
imaginaries if and only if the absolute Galois group of the relative
algebraic closure inside $F$ of the field generated by $C$ is isomorphic
to $\hat\zee$. It may therefore happen that ${\rm Th}(F)$ eliminates
imaginaries in the ring language, but it may also happen  that extra elements are needed,
for instance  if $F$ contains $\rat^{alg}$. The following lemma will therefore be
useful when dealing with $\aleph_\epsilon$-saturation.

\begin{lem}\vlabel{lem0} Let $F$ be an $\aleph_\epsilon$-saturated pseudo-finite
  field and $a$ a finite tuple in $F$. Then there is a finitely
  generated subfield $A$ of $F$ containing $a$ and such
  that $G(A^{alg}\cap F)\simeq \hat\zee$. 
\end{lem}

\prf Let $k$ be the relative algebraic closure inside $F$ of the subfield 
generated by $a$, and consider $k(t)$, where $t$ is transcendental over
$k$. Let $Q_0$ be the set consisting of all integers $n$ which are either
prime numbers or $4$ and such that $G(k)$ does not have a quotient isomorphic
to $\zee/n\zee$. If char$(k)\neq 0$, we let $Q=Q_0\setminus \{4\}$, and
if char$(k)=0$, we let $Q=Q_0\setminus \{2\}$. If $Q$ is empty, then
$G(k)\simeq \hat\zee$, and we are done. So, we assume that $Q$ is
non-empty.\\[0.05in]  
By Proposition 16.3.5 of \cite{FJ}, for each $n$, $k(t)$ has a Galois
extension $L_n$ which is regular over $k$ and with
$\gal(L_n/k(t))=\zee/n\zee$. Let $L$ be the field composite of all
$L_n$, $n\in Q$. Then $\gal(L/k(t))\simeq \prod_{n\in Q}\zee/n\zee$. 
Observe that $L\cap
k^{alg}=k$, because all $L_n$'s are regular extensions of $k$, and
Galois over $k(t)$ of
relatively prime order.  \\
Take a topological generator $\si_0$ of $\gal(L/k(t))$, and a topological
generator $\si_1$ of $G(k)$. Let
$\si\in G(k(t))$ extend $(\si_0,\si_1)\in \gal(Lk^{alg}/k(t))$ ($\simeq
\gal(L/k(t))\times G(k)$);  then the subfield $A$ of $k(t)^{alg}$ 
fixed by $\si$ is a regular extension of $k$, with Galois group
isomorphic to $\hat\zee$, since its Galois group is procyclic, projects onto $G(k)$, onto all 
$\zee/p\zee$ with $p$ a  prime,  and onto $\zee/4\zee$ if
char$(k)=0$. \\
By general properties of pseudo-finite fields and by
$\aleph_\epsilon$-saturation of $F$, there is a $k$-embedding $\varphi$
of $A$ inside  $F$, in such a way that
$\varphi(A)^{alg}\cap F=A$. This is classical, and follows from instance
from Lemma 20.2.2 in \cite{FJ}.


\begin{lem}\vlabel{lem1} Let $\kappa$ be an uncountable cardinal or
  $\aleph_\epsilon$, let
$K$ be a difference field, with $\fix(\si)(K)$ pseudo-finite and
$\kappa$-saturated. Then there is a model
 $\calu$ of ACFA containing $K$, which is $\kappa$-saturated, and with $\fix(\si)(\calu)=\fix(\si)(K)$.
\end{lem}

\prf (Compare with Afschordel's result \cite{A}). Let $\calu_1$ be a
$\kappa$-saturated model of ACFA containing $K$, and let
$\calu\subseteq \calu_1$ be maximal such that
$F:=\fix(\si)(\calu)=\fix(\si)(K)$. We will show that $\calu$ satisfies our
conclusion. First observe that $\calu$ is algebraically closed. 
Let $A=\acl(A)\subset \calu$ be small, let
$p\in S_1(A)$. Then $p$ is realised in $\calu_1$, and we take some $a\in
\calu_1$ realising $p$, with $\SU(a/\calu)$ minimal. Let $B\supset A$ be
small and 
such that $a\dnfo_B\calu$ and replace $p$ by $tp(a/B)$. \\[0.05in]
If $tp(a/\calu)\perp^a \fix(\si)$, then $\calu(a)_\si^{alg}$ has the same
fixed field as $\calu$: indeed, $\calu(a)_\si^{alg}$ and
$\fix(\si)(\calu_1)$ are linearly disjoint over their intersection,
which is contained in $\calu$ and therefore in $K$.  So by maximality of $\calu$,
$a\in\calu$. \\[0.05in]
Assume now that $tp(a/\calu)\not\perp^a \fix(\si)$. Then there is some
small $C\subset \calu$ containing $B$, and a realisation $a'$ of
$tp(a/B)$ such that $C(a')_\si\cap \fix(\calu_1)$ contains some element
$b$ not in $\calu$. We may and will assume that $\fix(\si)(C)$ has
absolute Galois group isomorphic to $\hat\zee$ (by Lemma \ref{lem0}). But as $F$ is $\kappa$-saturated, $tp_{F}(b/C\cap F)$ is realised in $F$, by some
$b_1$. Then $b_1$ realises $tp(b/C)$ (see the first paragraph of \ref{fix2}). Thus, by
$\kappa$-saturation of $\calu_1$, there is some $a_1\in\calu_1$ such that
$tp(a_1,b_1/C)=tp(a',b/C)$. But then $a_1$ realises $p$, and
$\SU(a_1/\calu)\leq \SU(a'/B)-\SU(b/C)<\SU(a/B)$, which gives us the desired
contradiction. \\[0.05in]
So in both cases, $p$ is realised in $\calu$. 

\begin{cor}\vlabel{corlem1} Let $\kappa$ be as above, and $K$ an
  algebraically closed difference field with $Fix(\si)(K)$
  $\kappa$-saturated.  If  $\calu$ is a
  $\kappa$-prime model of ACFA over $K$ then
  $\fix(\si)(\calu)=\fix(\si)(K)$. 
\end{cor}

\begin{lem}\vlabel{lem2} Let 
  $\calu$ be an $\aleph_\epsilon$-saturated
  model of ACFA of characteristic $0$, and let $K$ be an algebraically closed difference
  subfield of $\calu$ which contains $F:=\fix(\si)(\calu)$. Let
  $p=tp(a/K)$ be qf-internal to $F$, $p\perp^a\fix(\si)$, and assume
  $\si(a)\in K(a)$. Then there are  
  a (very) small $A\subseteq K$ 
  and a tuple $b$ of realisations of $p$ such
  that:
\begin{enumerate}
\item
  $FA(b)$ contains all realisations (in $\calu$) of
  $qftp(a/A)[\ell]$, for any $\ell\geq 1$.
\item If $b'$ realises $qftp(b/A)[m]$ for some $m\geq 1$, then
  $FA(b')$ contains all realisations  of $qftp(a/A)[\ell]$ for
  $\ell\geq 1$. 
\item $tp(a/A)\vdash tp(a/K)$, and $tp(b/A)\vdash tp(b/K)$.

\end{enumerate}
\end{lem}

\prf Let $k\subset K$ be small such that $a\dnfo_kK$ and
$\gal(\fix(\si)(k)^{alg}/\fix(\si)(k))$ is isomorphic to
$\hat\zee$. Then $\si(a)\in k(a)$ and $kF$ contains $\fix(\si^\ell)(\calu)$ for all
$\ell\geq 1$. By
assumption, there is some small  $B$ (in  $\calu$, by
$\aleph_\epsilon$-saturation of $\calu$) independent from $a$ over $k$, such
that $a\in BF$. Hence, there is a tuple $c$ in
$B(a)_\si\cap F=B(a)\cap F$ such that $B(a)=B(c)$.  Let $D=\Cb(a,c/B)$. Then
$D(c)=D(a)$, and  $D\subset
k(c_1,a_1,\ldots,c_n,a_n)$ for some independent realisations $(c_i,a_i)$ of
$qftp(c,a/B)$ (in some elementary extension of $\calu$). By $\aleph_\epsilon$-saturation
of $\calu$, we may assume that $(c_1,a_1,\ldots,c_n, a_n)$ is in
$\calu$, and is independent from $(c,a)$ over $D$. We let $b=(a_1,\ldots,a_n)$,
$A=\acb(k,c_1,a_1,\ldots,c_n,a_n/K)$; then $D\subset kF(b)$,
and $A$ is small. As $A$ contains  $c_1,\ldots,c_n$ ($\in F\subset K)$ and $k$, we also
have $D\subset A(b)$, whence $a\in FA(b)$. Note  
that $a\dnfo_kA$ since $A\subset K$.  \\[0.05in]
If $a'\in\calu$ realises $qftp(a/A(b))$, then the difference fields
$D(a)$ and $D(a')$ are isomorphic; hence there is some $c'\in D(a')\cap F$
such that $D(c')=D(a')$, i.e.: $a'\in FA(b)$. \\[0.05in]
Let $a'$ be an arbitrary realisation of $qftp(a/A)$, and let 
$b'$ be a realisation of $qftp(b/A)$, which is independent from
$(b,a')$ over $A$. By the previous paragraph (as $b'$ consists of $n$
realisations of $qftp(a/A(b))$) we know that $b'\in
FA(b)$. 
The difference fields $A(b)$ and $A(b')$ are $A$-isomorphic, and this
isomorphism extends to an isomorphism of difference fields
$A(b,a)\to A(b',a')$. Hence, $a'\in FA(b')\subseteq FA(b)$, as
desired. If $a'$ realises $qftp(a/A)[\ell]$ and is independent from $D$
over $k$, then the $\si^\ell$-difference fields $D(a')$ and $D(a)$ are
isomorphic over $D$. Let $f(x)$ be the tuple of rational functions over
$D$ such that $f(a)=c$; then $\si^\ell(f(a'))=f(a')$ and
$D(a')=D(f(a'))$. Hence $a'$ belongs to $FA(b)$. An argument similar
to the one given in the first case shows it for arbitrary realisation of
$qftp(a/A)[\ell]$ and shows (1). \\[0.1in]
Note that we have in fact shown that $FA(b')=FA(b)$, and so the
conclusion of (1) also holds for $b'$. An easy argument allows to remove
the assumption that $b'$ is independent from  $b$ over $A$: let $b''$
realise $qftp(b/A)$, independent from $(b,b')$ over $A$; then by the
proof of the first part: $FA(b'')=FA(b)$ and
$FA(b')=FA(b'')$. \\
Working in $\calu[\ell]$, and noting that if $m$ divides $\ell$, then
the realisations of $qftp(a/A)[m]$ also realise $qftp(a/A)[\ell]$,
part (1) gives (2). \\[0.1in]
%
For the proof of (3), we will first show that every realisation $b'$ of
$qftp(b/A)[\ell]$ (in $\calu$) is independent from $K$ over $A$. Indeed, by (2),
we know that $FA(b)=FA(b')$, and therefore $FK(b)=FK(b')$. This
implies that $\trdeg(b/K)=\trdeg(b'/K)$, and therefore that
$b'\dnfo_{A}K$. As $\calu$ is $\aleph_\epsilon$-saturated and $A$ is small, this shows that if $d\in K$, then
$$qftp(b/A)[\ell]\perp^a qftp(d/A)[\ell]. $$ 
By Proposition~4.9 of \cite{CH},
if $tp(b/A)\not\vdash tp(b/K)$, then there would be some tuple $d\in K$ and
integer $\ell\geq 1$ such
that $tp(d/A)[\ell]\not\perp^a tp(b/A)[\ell]$. But as we just saw,
this is impossible, and this gives us (3). (This is where the
characteristic $0$ assumption is crucial).

\begin{rem} In the above notation, note that if $\calu\prec\calu'$ and
  $F'=\fix(\si)(\calu')$, then $F'A(b)$ contains all realisations of
  $qftp(a/A)[\ell]$ in $\calu'$, for any $\ell>1$. 
  \end{rem} 

\begin{lem}\vlabel{lem4}
Let $K, A, b, \calu$ be as in Lemma \ref{lem2}, and let $L$ be a
difference subfield of $\calu$ containing $K$. Then there is a small
$A'$ containing $A$ such that $tp(b/A')\vdash tp(b/L)$. In
particular, $tp(a/A')\vdash tp(a/L)$. 
\end{lem}

\prf Let $A'\subset L$ be small, containing $A$ and such that
$b\dnfo_{A'}L$. Then the proof of (3) works. 

\begin{cor} \vlabel{corlem2} Let $K$ and $\calu$ be as in Lemma \ref{lem2}, and $p$ be a
  type which is almost internal to $\fix(\si)$. Then any $K$-indiscernible sequence $(a_i)$ of
  realisations of $p$ in $\calu$ is finite. 
\end{cor}

\prf Let $(a_i)_{i<\omega}$ be a sequence of realisations of $p$ in
$\calu$ which is $K$-indiscernible. 
Then either $a_0\in K$, or
$tp(a_0/K)$ is almost orthogonal to $\fix(\si)$ (since $K$ contains
$F:=\fix(\si)(\calu)$). By Lemma \ref{internality-3} there is $a'_0\in K(a_0)_\si$ such that
$\si(a'_0)\in K(a'_0)$, $a_0\in K(a'_0)^{alg}$ and $tp(a'_0/K)$ is qf-internal to
$\fix(\si)$. It suffices to show the result for $p=tp(a'_0/K)$. 
Let $b$ be the finite tuple of realisations of
$tp(a'_0/K)$ given by Lemma \ref{lem2}. If $n>d=\trdeg(K(b)/K)$
and $tp(a'_i,a_i/K)=tp(a'_0,a_0/K)$, then
we know that $a'_n\in K(a'_0,\ldots,a'_{d-1})^{alg}$ (because $K\supset
F$). Hence the sequence is finite.

\begin{defn} We call a type $p$ over a set $A$ {\em acceptable}
  (in $K\supset A$) if $A$ is the algebraic closure of a finite
  tuple, and 
  either $\SU(p)=1$ and $p$ is one-based, or $p$ is qf-internal to
  $\fix(\si)$,  almost orthogonal to $\fix(\si)$, and  if $b$ realises  $p$  then
  $\si(b)\in A(b)$,  $tp(b/A)\vdash tp(b/K)$,
and the
  set of realisations of $qftp(b/A)[\ell]$ for $\ell\geq 1$, in some model $\calu$ of
  ACFA containing $K$, is contained in $A(b)\fix(\si)(\calu)$. 
\end{defn}

\begin{notation} Let $p$ be a one-based type of SU-rank $1$ over the very small set
  $A$. If $A\subset B\subset C$, we denote by $p|B$ the unique non-forking
  extension of $p$ to $B$, and by $\dim_Bp(C)$ the cardinality of a
maximal  $B$-independent subset of realisations of $p|B$ in $C$. 
\end{notation}

\begin{lem}\vlabel{onebased1} Let $p$ be an acceptable
  one-based type over the very small $A$, let $K$ be an
  algebraically closed difference field containing $A$. We work in a
  sufficiently saturated model $\calu$ of ACFA. Let $\kappa$ be
  an infinite cardinal. 
  \begin{enumerate}
    \item If $K$ contains $\kappa$  many
  $A$-independent realisations of $p$, then the non-forking extension of $p$
  to $K$ is not $\kappa$-isolated, and conversely.
  \item One of the following holds:  \begin{itemize}
      \item[(a)] There is some $n<\omega$ and realisations
        $a_0,\ldots,a_{n-1}$ of  $p|K$ such that \hfill \break 
        $\dim_A\,p(\acl(Ka_0\ldots,a_{n-1}))\geq \kappa>\dim_A\,p(K)$. Moreover, if
        $n$ is 
        minimal with this property,
 then
        $tp(a_0,\ldots,a_{n-1}/K)$ is $\kappa$-isolated (but
 $p|\acl(K a_0\ldots,a_{n-1})$ is not).
 \item[(b)] Not case (a). If $B$ is a set of $K$-independent
          realisations of $p|K$ of size $\lambda<\kappa$, then
          $\dim_A \,p(\acl(KB))<\kappa$. 
      \end{itemize}
  
  \end{enumerate}
  \end{lem}

  \prf (1) If $C=\acl(C)\subset K$ is small, then $C$ contains 
  $<\kappa$ $A$-independent realisations of $p$, so that the non-forking extension
  of $p$ to $C$ will be realised in $K$, and $p|K$ is not $\kappa$-isolated. The converse is clear: the
  non-forking extension of $p$ to $K$ is implied by its restriction to
  $\acl(A,p(K))$. \\[0.05in]
  (2) (a) is clear by (1) and because $\dim$ is additive. So, assume that there is no such $n$, and let
  $B$ be as in (b), and  $(a_i)_{i<\lambda}\subset B$  a sequence of independent over $K$
realisations of $p$, and assume that
$\lambda< \dim_A p(\acl(KB))=\mu\geq \kappa$. 
So $\acl(K,a_i\mid i<\lambda)$
contains a set $C$ consisting of $\mu$ many $A$-independent realisations of $p$. Then for each $c\in C$, there is
some finite $I_c\subset \lambda$ 
 such that $c\in \acl(K a_i\mid i\in I_c)$. As $\lambda<\mu$, some
 set $I_c$ appears $\mu$ times. Thus $\dim_A\, p(\acl(Ka_i\mid i\in
 I_c))=\mu\geq \kappa$, which contradicts our
 assumption.

 \begin{rem}\vlabel{trans1} Let $p$ be the generic $1$-type over
   $K$, and $\kappa$ an infinite cardinal. Then $p$ is $\kappa$-isolated
   if and only if $\Delta(K)<\kappa$. This follows easily from the
   description and properties of the generic types, see \ref{generic}. 
   \end{rem}

\begin{defn} Let $K=\acl(K)\subset L=\acl(L) \subset \calu$. We say that $L$ is {\em normal over $K$} (in $\calu$) is
  whenever $a$ is a tuple in $L$, then $L$ contains all realisations of
  $tp(a/K)$ in $\calu$. 
\end{defn}

\begin{lem}\vlabel{lem5} Let $\kappa$ be an uncountable cardinal or
  $\aleph_\epsilon$, let $K\subseteq L$  be algebraically closed difference
  subfields of $\calu$, where $\calu$ is $\kappa$-saturated, and
  $\fix(\si)(\calu)\subset K$. Assume that $\calu$ is
  $\kappa$-atomic over $K$.
  \begin{enumerate}
 \item Let $B\subset \calu$ be transformally independent over $K$, and
   assume that either $|B|<\kappa$, or that $B$ is a tranformal
   transcendence basis of $\calu$ over $K$. Then $\calu$ is
   $\kappa$-atomic over $\acl(KB)$.    
\item If $L$ is normal over $K$ then
  $\calu$ is $\kappa$-atomic over $L$.
\end{enumerate}  
\end{lem} 

\prf (1) If $B\neq \emptyset$ there is nothing to prove, so suppose it is
not.  Then
$\Delta(K)<\kappa$ by \ref{trans1}. Let $a$ be a finite tuple in $\calu$, and let $b\subset B$ be a
finite tuple such that $a\dnfo_{Kb}B$. Let $c\subset a$ be a tranformal
transcendence basis of $K(a,b)_\si$ over $K(b)_\si$ (and therefore also
over $K(B)_\si$). If $c\neq \emptyset$, then $|B|<\kappa$, $
\Delta(K(B))<\kappa$, and therefore
$tp(c/\acl(KB))$ is $\kappa$-isolated. Moreover, as $a$ is transformally
algebraic over $K(b,c)_\si$, and $B$ is purely transformally
transcendental over $K(b,c)_\si$, $tp(a/\acl(Kbc))$ and
$tp(B/\acl(Kbc))$ are orthogonal, and by stationarity of
$tp(B/\acl(Kbc))$, we get that $tp(B/\acl(Kbc))\vdash
tp(B/\acl(Kba))$. By symmetry, 
$$tp(a/\acl(Kbc))\vdash tp(a/\acl(KBc)).$$  But $tp(a,b,c/K)$ is
$\kappa$-isolated, and this gives the result. \\[0.05in]
(2) Let $a$ be a finite tuple in $\calu$, and consider $tp(a/L)$. Let
$d\subset a$ be maximal transformally independent over $L$. If $d\neq
\emptyset$, then $d$
is  transformally independent over $K$, which implies that
$\Delta(L/K)=0$ (by normality of $L/K$), and that
$\Delta(K)=\Delta(L)<\kappa$  (by $\kappa$-isolation of
$tp(d/K)$). Therefore  $tp(d/L)$ is
$\kappa$-isolated. \\
If $\Delta(L/K)\neq 0$, note that by
normality of $L$ over $K$ in $\calu$, every element of the tuple $a$ which is not
in $L$ is transformally algebraic over $K$. So, replacing $a$ by
$a\setminus L$, we may assume they are all transformally algebraic over
$K$, i.e., that $\SU(a/K)<\omega$. We then let $d=\emptyset$. \\[0.05in]
In both cases, by  (1), 
$\calu$ is $\kappa$-atomic over $\acl(Kd)$, and the normality of $L$
over $K$ implies the normality of 
$\acl(Ld)$ over $\acl(Kd)$. Working over $\acl(Kd)$, we may
therefore assume that $a$ and $D:=\acb(a/L)$ are 
 transformally algebraic over $K$. \\
We use induction on $\SU(a/L)$, and using the
semi-minimal analysis, we find $b\in \acl(Da)$ such that
$tp(a/\acl(Db))$ is either one-based of SU-rank $1$, or almost-internal
to $\fix(\si)$. \\
If
$tp(a/\acl(Db))$ is almost-internal to $\fix(\si)$, then so is
$tp(a/\acl(Lb))$. By Lemma \ref{internality-3}, there is $a'\in\acl(Lba)$
such that $a\in\acl(Lba')$ and $tp(a'/\acl(Lb))$ is qf-internal to
$\fix(\si)$. By Lemma \ref{lem2}, there is a very small $D'\supseteq D$
such that $tp(a'/\acl(D'b))\vdash tp(a'/\acl(Lb))$, and we may choose it so
that $a\in \acl(D'ba')$. This shows  that $tp(a/\acl(Lb))$ is
$\kappa$-isolated, and therefore so  is $tp(a/L)$.  \\
So assume  that $p:=tp(a/\acl(Db))$ is one-based of $\SU$-rank $1$, and
let $c$ be a tuple containing $b$ such that
$\acl(Db)=\acl(c)=:C$.  We need to show that
$\dim_C\, p(\acl(Lc))<\kappa$. As  $\calu$ is
$\kappa$-atomic over $K$, we know that $tp(a,c/K)$ is
$\kappa$-isolated, and therefore $\dim_C\, p(\acl(Kc))<\kappa$. 
So,
if $\dim_C\, p(\acl(Lc))\geq \kappa$, then there is some $a'\in
\acl(Lc)\setminus \acl(Kc)$ realising $p$. Recall that by our earlier step, $c,a'$ are transformally
algebraic over $K$, and therefore so is $e=\Cb(c,a'/L)$.  Consider now $\acl(Kca')\cap
\acl(Ke)=:E\subset L$;  by Proposition~3.1 
of \cite{ACFAmod}, $tp(e/E)$ is almost-internal to $\fix(\si)$, and
therefore orthogonal to all one-based types. As
$tp(a'/Kc)$ is one-based, and $a'\in \acl(Kce)\setminus \acl(Kc)$, it
follows that $e\in E$, since internality to $\fix(\si)$ and
non-orthogonality to a one-based type imply algebraicity. That is, $e\in
\acl(Kca')\cap L$, and as $a'\notin \acl(Kc)$, the tuples $a'$ and $e$ 
are equi-algebraic over $\acl(Kc)$. Hence 
$\acl(Kca)$ contains a
realisation of $tp(e/\acl(Kc))$, because $tp(a/\acl(Kc))=tp(a'/\acl(Kc))$.  But this contradicts the normality of
$\acl(Lc)$ over $\acl(Kc)$.
So, $\dim_C(p(\acl(Lc))<\kappa$, and $tp(a/Lb)$ is $\kappa$-isolated.

\begin{thm}\vlabel{thm1} Let $\kappa$ be an uncountable cardinal or
  $\aleph_\epsilon$, and let $K$ be an algebraically closed difference field of
characteristic $0$ such that $F:=\fix(\si)(K)$ is pseudo-finite and is
$\kappa$-saturated. 
\begin{enumerate}
\item Then there is a $\kappa$-prime model $\calu$ over $K$. 
\item Furthermore, $\calu$ is $\kappa$-atomic
  over $K$, and every
sequence of $K$-indiscernibles has length $\leq \kappa$ (i.e., if
$\kappa=\aleph_\epsilon$, $\leq \aleph_0$; by convention, if $\kappa$ is
meant as a cardinal, then $\aleph_\epsilon$ will mean $\aleph_0$).
\end{enumerate}
\end{thm}

\prf 
By Lemma \ref{lem1}, there is a $\kappa$-saturated model $\calu_1$ of ACFA containing $K$ and with fixed field
$F=\fix(\si)(K)$.  We will construct a submodel $\calu$ of $\calu_1$, which
is $\kappa$-prime over $K$ and satisfies
(2). \\[0.05in] 
{\bf Step 0}. Taking care of the transformal transcendence degree. \\
If the transformal
transcendence degree of $K$ is $<\kappa$, then as any
$\kappa$-saturated model of ACFA has
transformal transcendence degree at least $\kappa$, we enlarge $K$
as follows: let $B\subset \calu_1$ be a set which is transformally
independent over $K$ and of cardinality $\kappa$; by 2.11 of \cite{CH}, this condition completely
determines the $K$-isomorphism type of $K(B)_\si^{alg}$, and therefore
any $\kappa$-prime model will contain a
$K$-isomorphic copy of $K(B)_\si^{alg}$. We let $K_0=K(B)_\si^{alg}$. We
need to show (2). Each finite subset of $B$ realises a
$\kappa$-isolated type over $K$, since the transformal
transcendence degree of $K$ is $<\kappa$. Moreover, every tuple in $K_0$ realises an isolated type
over $K(B)_\si$; hence $K_0$ is $\kappa$-atomic over $K$. It is also
$\kappa$-constructed over $K$.\\
Let
$(a_i)_{i<\lambda}\subset K_0$ be  a $K$-indiscernible sequence, $\lambda$ a cardinal. If the $a_i$'s are transformally independent
over $K$, then we know that $|\lambda|\leq \kappa$. If not,
then by indiscernibility, the transformal transcendence degree of $K(a_i\mid
i<\lambda)_\si$  over $K$ is finite, and we choose a finite subset $c$
of $B$ such that $K(a_i\mid
i<\lambda)_\si$ is transformally algebraic over $K(c)_\si$. As the
elements of $B$ are transformally independent over $K$, this implies that all
$a_i$'s  are in fact algebraic over $K(c)_\si$. Consider now
$D:=\acb(c/Ka_i\mid i<\lambda)$. For every $i$, we know that $a_i\in
K(c)_\si^{alg}$, and therefore by definition of $D$, $a_i\in
D(c)_\si^{alg}$. But $c$ is finite, $D$ is contained in the algebraic
closure of a finite set (by \ref{cb}), and therefore  $D(c)_\si^{alg}$
is countable. Hence so is $\lambda$. 
 This shows condition (2) for the extension
$K_0/K$. \\[0.1in]
We will build a
sequence $K_n$, $n<\omega$,   of algebraically closed
difference subfields of $\calu_1$, such that:\\
(i) if $p$ is an acceptable
type over a very small $A\subset K_n$, then  $K_{n+1}$ contains
$\kappa$-many $A$-independent realisations of $p$;\\
(ii) $K_{n+1}$ is $\kappa$-constructed over $K_n$. \\[0.05in]
We let $K_0=K$ if the transformal transcendence degree of $K$ is $\geq
\kappa$, and $K(B)_\si^{alg}$ as in step $0$ otherwise. We
assume $K_n$ constructed, we wish to build $K_{n+1}$. Let
$p_\beta$, $\beta<\lambda$, be an enumeration of the acceptable types in 
$K_n$, with corresponding very small bases $A_\beta$. \\[0.05in]
{\bf Step 1}. Defining $K_{n+1}=\bigcup_{\beta<\lambda}K'_\beta$\\
We  build the sequence $K'_\beta$ by
induction on $\beta$, and let $K'_0=K_n$. If $\beta$ is a limit
ordinal, then we let $K'_\beta=\bigcup_{\gamma<\beta}K'_\gamma$, and $K_{n+1}=K'_\lambda$. We will
build them so that $K'_{\beta+1}$ satisfies the following:\\
(i') $K'_{\beta+1}$ contains $\kappa$-many $A_\beta$-independent
realisations of $p_\beta$;\\
(ii') $K'_{\beta+1}$ is $\kappa$-constructed over $K'_\beta$.\\[0.05in]
Assume
$K'_\beta$ constructed. 
If $p_\beta$ has $\kappa$-many $A_\beta$-independent
realisations in $K'_\beta$, then we let
$K'_{\beta+1}=K'_{\beta}$. Otherwise, we need to distinguish the two
cases:\\[0.05in]
{\bf Case 1}.  $p_\beta$ is one-based.\\
Let $a_i$, $i<\kappa$, be a sequence of $K'_\beta$-independent
realisations of $p_\beta$ (a priori, in some elementary extension of
$\calu_1$). By Lemma \ref{onebased1}, either there is
$n<\omega$ such that $\acl(K'_\beta,a_i\mid i<n)$ contains $\kappa$-many
$A_\beta$-independent realisations of $p_\beta$; in that case, taking a minimal such
$n$,  $tp(a_0,\ldots,a_{n-1}/K'_\beta)$ is $\kappa$-isolated and
therefore realised in $\calu_1$, so that we may assume
$a_0,\ldots,a_{n-1}\in\calu_1$ and we set 
$K'_{\beta+1}=\acl(K'_\beta,a_i\mid i<n)$. Then (i') and (ii') follow. \\
If there is no such $n$, by the same reasoning we may assume the $a_i$'s
are in $\calu_1$ and we define $K'_{\beta+1}=\acl(Ka_i\mid
i<\kappa)$. Then (i') and (ii') again are satisfied.   \\[0.05in]
{\bf Case 2}. Not case 1. \\
Let $a_\beta\in \calu_1$ realise $p_\beta$,
$K'_{\beta+1}=K'_{\beta}(a_\beta)^{alg}$. By assumption on $p_\beta$,
$tp(a_\beta/A_\beta)\vdash tp(a_\beta/K_{n})$. By
Lemma~\ref{lem4}, there is a very small subset $B$ of $K'_\beta$ which
contains 
$A_\beta$ and is such that $tp(a_\beta/B)\vdash tp(a/K'_\beta)$. So,
$tp(a_\beta/K'_\beta)$ is $\kappa$-isolated. We let
$K'_{\beta+1}=K'_\beta(a_\beta)_\si^{alg}$. We know that
$FK'_\beta(a_\beta)_\si$ contains all realisations of
$tp(a_\beta/B)$ in $\calu_1$. But $\calu_1$ is $\kappa$-saturated, it
therefore contains $\kappa$ independent realisations of
$tp(a_\beta/A_\beta)$, which shows (i'). \\[0.05in]
We now define  $\calu=\bigcup_{n\in \omega}K_n$. \\[0.1in]
{\bf Step 2}. Show that $\calu$ is $\kappa$-saturated.\\
Let $C\subset \calu$ be small, and $p$ a $1$-type over $C$, realised by
$a$ in $\calu_1$. If $\SU(p)=\omega$, then $a$ is transformally
transcendental over $C$; as $C$ is small, $K_0$ will contain a
realisation of $p$. So we may assume that $\SU(p)<\omega$, and the proof is by
induction on $\SU(p)$: we assume that for any small $D$, any $1$-type $q$
over $D$ of smaller SU-rank than $p$ is realised in $\calu$. \\
If
$\SU(p)=0$ there is nothing to prove, as $p$ is realised in $C$. If there
is some $b\in C(a)_\si^{alg}$ such that $0<\SU(b/C)<\SU(a/C)$, then we get
the result 
 by induction: $tp(b/C)$ is realised by some $b'\in\calu$, and
there is $a'\in\calu$ such that $tp(a',b'/C)=tp(a,b/C)$, since
$\acl(Cb')$ is small and $\SU(a/Cb)<\SU(p)$. \\[0.05in]
Hence we may assume that there is no such $b$, whence $p$ is either one-based of SU-rank $1$, or almost-internal to
$\fix(\si)$ (by the semi-minimal analysis \ref{semi-min}). We need to
distinguish  three cases. \\
{\bf Case 1}. $p$ is one-based of SU-rank $1$. \\ 
Let $A\subset C$ be very small such that $p$ does not fork over $A$. Let
$n<\omega$ be such that $A\subset K_n$; then $p$, being acceptable,
occurs as a $p_\beta$, and is therefore realised in $K_{n+1}$. \\[0.05in]
{\bf Case 2}. $p$ is realised in $\fix(\si)$. \\
If $a\in\fix(\si)$, we saw in \ref{fix2} that $tp_F(a/C\cap F)\vdash
tp(a/C)$. The saturation hypothesis on $F$ then gives the result: $p$ is
realised in $F$. \\
{\bf Case 3}. Assume now that $p\perp^a \fix(\si)$, $p$ almost-internal
to $\fix(\si)$. \\
By Lemma \ref{internality-3}, there is $a_1\in C(a)_\si$ such that $tp(a_1/C)$ is qf-internal to
$\fix(\si)$, $\si(a_1)\in C(a_1)$, and $a\in C(a_1)^{alg}$. We may replace $p$ by $tp(a_1/C)$,
i.e., assume that $p$ is qf-internal to $\fix(\si)$. Let $C_0\subset C$
be very small such that $p$ does not fork over $C_0$.  By Lemma \ref{lem2}
there is a tuple $b$ of realisations of $p$ and a very  small $D$
containing $C_0$, contained in $\acl(CF)$, such that $FD(b)$ contains all
realisations of $qftp(a/D)$, and $tp(b/D)\vdash tp(b/\acl(CF))$. Thus,
$tp(b/D)$ is acceptable, and if $n$ is such that $D\subset K_n$, then
$p$ in realised in $K_{n+1}$.\\[0.1in]
{\bf Step 3}. $\calu$ is $\kappa$-prime over $K$.\\
This is clear, by Remarks \ref{atom}(2)(3). \\[0.1in]
%
%
{\bf Step 4}.  $\calu$ is $\kappa$-atomic 
over $K$. \\ 
When $\kappa$ is regular or $\aleph_\epsilon$, then this folllows from
$\calu$ being $\kappa$-constructed over $K$.  
The proof in the singular case is a little more delicate, and is done by
induction. We already saw that $K_0$ is $\kappa$-atomic over $K$. Let
$a$ be a finite tuple in $\calu$, and (in the notation of Step 1), choose $n$ minimal 
such that $a\in K_{n+1}$, and $\beta$ minimal such that $a\in
K'_{\beta+1}$. If $n=-1$, there is nothing to prove (by Step 0), so assume $n\geq
0$. By definition of $K'_{\beta+1}$, there are a tuple $b$  in $K'_\beta$
and a tuple $c$ of realisations of $p_\beta$ such that $a\in
\acl(Kbc)$. We may assume that $\acl(Kb)$ contains $A_\beta$, and that
$c\dnfo_{Kb}K'_\beta$. By induction hypothesis, $tp(b/K)$ is
$\kappa$-isolated, and it therefore suffices to show that
$tp(c/\acl(Kb))$ is $\kappa$-isolated (by \ref{atom}(1)). If $p_\beta$ is
qf-internal to $\fix(\si)$  then we know by Lemma~\ref{lem2} that there is some very
small $D\subset \acl(Kb)$ such that $tp(c/D)\vdash
tp(c/\acl(Kb))$, and we are done. \\
If $p_\beta$ is one-based, then we may
assume that the elements of the tuple $c$ are independent over
$K'_\beta$, maybe at the cost of increasing $b\in K'_\beta$. Then, by the construction of $K'_{\beta+1}$ in Step 1, we know
that $tp(c/K'_\beta)$ is $\kappa$-isolated, so that 
if $c'$ is a proper subtuple of $c$ (consisting of realisations of $p_\beta$), then  $\dim_{A_\beta}\,p_\beta(\acl(K'_\beta
c'))<\kappa$. In particular, $\dim_{A_\beta}\,p(\acl(Kbc'))<\kappa$, and
$tp(c/\acl(Kb))$ is $\kappa$-isolated (by Lemma
\ref{onebased1}). \\[0.05in]
{\em Remarks}. (Notation as in step 1 and above) The same proof shows
that $\calu$ is $\kappa$-atomic over 
each $K_n$, and over each $K'_\beta$. \\
Moreover, the fact that $\calu$ is $\kappa$-atomic over $K'_\beta$
implies that $p_\beta(\calu)\subset K'_{\beta+1}$. \\[0.1in] 
{\bf Step 5}. If $(b_i)_{i<\lambda}\subset \calu$ is $K$-indiscernible,
with $\lambda$ a cardinal, 
then $\lambda\leq \kappa$. \\
By supersimplicity,  for some $n<\omega$ the elements $b_i, n<i<\lambda$, are
independent over $K(b_0,\ldots,b_n)_\si$. If
$\SU(b_{n+1}/Kb_0,\ldots,b_n)\geq \omega$, then the tuple $b_{n+1}$
contains an element which is transformally transcendental over $K$, and as the transformal
transcendence degree of $\calu$ over $K$ is $\leq \kappa$, we get
$\lambda\leq \kappa$. So we may assume $\SU(b_{n+1}/\acl(Lb_0\ldots,b_{n}))<\omega$.\\ 
Let $L=\acl(Kb_0,\ldots,b_n)$. Then the sequence $(b_i)_{n<i<\lambda}$
is indiscernible over $L$. Note that the sequence $\acl(Lb_i)$,
$n<i<\lambda$, is also indiscernible over $L$ under a suitable
enumeration of each $\acl(Lb_i)$. Hence, if $c_{n+1}\in \acl(Lb_{n+1})$,
there are $c_i\in \acl(Lb_i)$, 
$n+1<i<\lambda$, such that the sequence $(c_i)_{n<i<\lambda}$ is
indiscernible over $L$. 
Using the semi-minimal analysis \ref{semi-min} we may therefore assume that either $tp(c_i/L)$ is one-based of
SU-rank $1$, or that $tp(c_i/L)$ is almost-internal to $\fix(\si)$. If  $tp(c_i/L)$ is almost-internal to $\fix(\si)$, then the result follows by
Corollary \ref{corlem2}. The one-based case is a little more
complicated. \\[0.05in]
Towards a contradiction, assume that $\lambda>\kappa$,  that
$tp(c_{n+1}/L)$ is one-based of SU-rank $1$, and let
$C\subset L$ be a very small set   such that $tp(c_{n+1}/L)$ does not
fork over $C$, and set $p=tp(c_{n+1}/C)$. Then the tuples $c_i$, $n<i<\lambda$, form a Morley sequence over $C$
and over $L$. Let $N$ be $\kappa$-prime over
$M:=\acl(L,c_i \mid n<i<\kappa)$. We may assume that $N\prec
\calu$.\\[0.05in] 
     {\bf Claim}. $\calu$ is $\kappa$-prime over $L$. \\
     It suffices to show that $\calu$ is $\kappa$-constructed over $L$. To do
     that it is enough to show that each $LK_m$ is $\kappa$-constructed over
     $LK_{m-1}$.\\
If $m=0$ and $K_0\neq K$,  let $B_0$ be a finite subset of $B$  (the
transformal transcendence basis of $\calu$ over $K$)  such that
$b:=(b_0,\ldots, b_n)$ is independent from $K_0$ over $\acl(KB_0)$. In
particular, $b$ is transformally algebraic over $\acl(KB_0)$, and
therefore  $tp(B/\acl(KB_0)) \vdash tp(B/\acl(LB_0))$ (reason as in the
proof of Lemma \ref{lem5}(1)), and as
$tp(B_0/L)$ is $\kappa$-isolated, it follows that $K_0$ is $\kappa$-constructed
over $L$.\\
Assume now $n\geq 0$, and that we have shown that $K'_\beta$ is $\kappa$-constructed
over $L$. If $p_\beta$ is not one-based, then by Lemma \ref{lem4},
$tp(a_\beta/\acl(K'_\beta L))$ is $\kappa$-isolated, and we are
done. Assume now that $p_\beta$ is one-based; by construction there
is a set $(a_\alpha)_{\alpha<\mu}$ of $K'_\beta$-independent realisations of
$p_\beta|K'_\beta$ such that
$K'_{\beta+1}=\acl(K'_\beta,a_\alpha,\alpha<\mu)$, and either $\mu\in
\omega$ or $\mu=\kappa$. \\
If $\mu\in \omega$,  as $\calu$ is $\kappa$-atomic over $K'_\beta$, we
get 
$tp(a_0,\ldots,a_{\mu-1},b/K'_\beta)$ is $\kappa$-isolated and therefore
$LK'_{\beta+1}$ is $\kappa$-constructed over $LK'_\beta$. If $\mu=\kappa$, then
$\dim_{K'_\beta}p_\beta(\acl(
K'_\beta,b,a_\gamma\mid \gamma<\alpha))<\kappa$ for each $\alpha<\kappa$,  so that $K'_{\beta+1}$
is $\kappa$-constructed over $LK'_\beta$ (here we use that 
$p_\beta(\calu)\subset K'_{\beta+1}$ and that $b$ is finite). \\[0.05in]
Hence, $\calu$ being $\kappa$-prime
 over $L$, there is an $L$-embedding $f$ of $\calu$ into
$N$. So we have $L\subset f(\calu)\prec N\prec \calu$. 
As $\lambda>\kappa$ and the $c_i$'s are independent over $L$, there is some $n<j<\lambda$ such that
$f(c_j)\notin M$. But $\dim_M(p)\geq \kappa$, and by Lemma 
\ref{onebased1}, $p|M$ is not isolated. But $N$ is $\kappa$-atomic over $M$, and
$f(c_j)$ realises $p$ and is not in $M$, which gives us the desired
contradiction. 
This finishes the proof of (2) and of the theorem.


\begin{prop}\vlabel{lem6} Let $\kappa$ be an uncountable cardinal or
  $\aleph_\epsilon$, let $\calu$ and $\calu'$ be
  $\kappa$-saturated models of ACFA of characteristic $0$. Assume that
$\calu$ (resp. $\calu'$) contains an algebraically closed difference
field $K$ (resp. $K'$), over which it is $\kappa$-atomic and over which every
sequence of indiscernibles has length $\leq\kappa$. Assume moreover that 
  $F:=\fix(\si)(K)=\fix(\si)(\calu)$, $\fix(\si)(K')=\fix(\si)(\calu')$,
and that we have an isomorphism $f:K\to K'$. Let  $p$ be an acceptable type
over some very small $A\subset K$, and $p'=f(p)$. If $L=\acl(K p(\calu))$ and
$L'=\acl(K p'(\calu'))$, then $f$ extends to an isomorphism between $L$ and $L'$. \end{prop}

\prf Note that $p'$ is also acceptable, with very small basis $A'=f(A)$.
If $p$ is not one-based, then this is clear  by Lemma \ref{lem2}: $L=\acl(Kb)$,
$L'=\acl(K'b')$ for some tuples  $b$ realising $p$ and $b'$ realising
$p'$. 
We extend $f|A$ to an isomorphism $g_0:\acl(Ab)\to
\acl(A'b')$ which sends $b$ to $b'$; as $tp(b/A)\vdash tp(b/K)$ and
$tp(b'/A')\vdash tp(b'/K')$, $g_0\cup f$ extends to an isomorphism
$\acl(Kb)\to \acl(K'b')$. \\
Assume now that 
$p$ is one-based. Any
$\kappa$-saturated model of ACFA containing $A$ will contain
(at least) $\kappa$ realisations of $p$ which are independent over $A$; hence so
do $\calu$ and $\calu'$. Let
$(a_i)_{i<\lambda}\subset \calu$ be a set of realisations of $p$ which
is maximal independent over $K$, with $\lambda$ a cardinal, and let
$(a'_i)_{i<\mu}\subset \calu'$ be defined analogously over $K'$. By Lemma \ref{onebased1} and our
hypothesis on the length of $K$-indiscernible sequences, either
$\lambda$ is finite, or $\lambda=\kappa$. If $\lambda=n<\omega$, then as
$tp(a'_0,\ldots,a'_{n-1}/K')=f(tp(a_0,\ldots,a_{n-1}/K))$,  it follows that $\acl(K'a'_0,\ldots,a'_{n-1})$ contains
$\kappa$-many independent realisations of $f(p)$, so that  $\mu\leq n$. The
symmetric argument gives $\mu=\lambda$. Define $g$ on $K(a_i\mid
i<\lambda)_\si$ by $g(a_i)=a'_i$, and extend to $L=\acl(Ka_i\mid
i<\lambda)$.

\begin{thm}\vlabel{thm2} Let $\kappa$ be an uncountable cardinal or
  $\aleph_\epsilon$, let $\calu$ and $\calu'$ be
  $\kappa$-saturated models of ACFA of characteristic $0$ containing an
  algebraically closed difference field $K$, with
  $F:=\fix(\si)(K)=\fix(\si)(\calu)=\fix(\si)(\calu')$. Assume that
  $\calu$ and $\calu'$ are $\kappa$-atomic over $K$, and that
  any sequence of $K$-indiscernibles in $\calu$ or in $\calu'$ has
  length $\leq \kappa$. Then $\calu\simeq _K\calu'$.
\end{thm}

\prf We start with the generic type: if the transformal transcendence
degree of $K$ is $\geq \kappa$, then $\calu$ and $\calu'$ are
transformally algebraic over $K$. If not, then let $D$ be a transformal
transcendence basis of $\calu$ over $K$, $D'$ a transformal
transcendence basis of $\calu'$ over $K$. They have the same
cardinality $\kappa$, and there is a $K$-isomorphism $K(D)_\si^{alg}\to
K(D')_\si^{alg}$. By Lemma~\ref{lem5}, $\calu$ and $\calu'$ still
satisfy the hypotheses over $K(D)_\si^{alg}$ and $K(D')_\si^{alg}$. Hence
 we may assume that both $\calu$ and
 $\calu'$ are 
transformally algebraic over $K$. We define by induction on $n$ an
increasing 
sequence $K_n$ of algebraically closed subfields of $\calu$ such that
for each $n$, if $p$ is an acceptable type over some (very small) $A\subset K_{n-1}$,
then $K_n$ contains all realisations of $p$ in $\calu$, and furthermore,
$K_n=\acl(K_{n-1}P)$ for the set $P$ of all realisations (in $\calu$) of acceptable types
 over some subset of $K_{n-1}$. Then each $K_n$ is normal over
 $K_{n-1}$ (and in fact over $K$), and so by Lemma \ref{lem5}, $\calu$ satisfies the hypotheses
 over $K_n$. Note also that $\calu=\bigcup_{n<\omega}K_n$. We let
 $L_n\subset \calu'$ be defined analogously. It then suffices to build
 a sequence $g_n$ of $K$-isomorphisms $K_n\to L_n$. \\
Assume $g_{n-1}$
 already built. Let $p_\beta$, $\beta<\lambda$, be an enumeration of all
 acceptable types over a subset of $K_{n-1}$, with associated small
 basis $A_\beta$. Note that $f(p_\beta), \beta<\lambda$, enumerates all acceptable types
 over subsets of $L_{n-1}$, since if $q$ is an acceptable type over the
very small 
 $C\subset L_{n-1}$, so is $g_{n-1}\inv(q)$ (over $g_{n-1}\inv(C)\subset K_{n-1}$). We
 build by induction on $\beta<\lambda$ an increasing sequence $K'_\beta$
 of algebraically closed difference subfields of $\calu$ such that
 $K'_\beta$ contains all realisations in $\calu$ of $p_\gamma$ for all
 $\gamma<\beta$. Assume we have extended $g_{n-1}$ to an isomorphism
 $f_\beta:K'_\beta\to L'_\beta$, where $L'_\beta$ contains all
 realisations in $\calu'$ of $g_{n-1}(p_\gamma)$ for all $\gamma<\beta$. As
 $\calu$ is $\kappa$-atomic over $K_{n-1}$, it is also $\kappa$-atomic 
 over $K'_\beta$ (by Lemma \ref{lem5}), and similarly, $\calu'$ is
 $\kappa$-atomic over $L'_\beta=f_\beta(K'_\beta)$. Extending $f_\beta$
 to an isomorphism $f_{\beta+1}:K'_{\beta+1}\to L'_{\beta+1}$ is given
 by Lemma \ref{lem6}. 

As remarked before, if $q$ is an acceptable type over some $A'\subset
L'_{n-1}$, then $g_{n-1}\inv(q)=p_\beta$ for some $\beta<\lambda$, and so
$L'_n$ contains $q(\calu')$, and $K'_n$ contains
$g_{n-1}\inv(q)(\calu)$. This finishes the induction step. Then
$g=\bigcup_{n<\omega}g_n$ is a $K$-isomorphism between $\calu$ and
$\calu'$.

\begin{thm} \vlabel{corthm}
  Let $\kappa$ be an uncountable cardinal or
  $\aleph_\epsilon$, let $K$ be an algebraically closed difference
  field of characteristic $0$,  with $\fix(\si)(K)$ pseudo-finite and
  $\kappa$-saturated. Then ACFA has a $\kappa$-prime model over $K$, and
  it is unique up to $K$-isomorphism.
  
\end{thm}

\prf This follows immediately from Theorem \ref{thm2}, as the properties
are preserved by elementary substructures.

\begin{rem} Note that the result also holds under the weaker hypothesis:
  $K$ algebraically closed, 
  $|\fix(\si)(K)|<\kappa$, and $\kappa^{<\kappa}=\kappa\geq \aleph_1$, so that the
  theory of pseudo-finite fields has a unique (up to $K$-isomorphism)
 saturated model of cardinality $\kappa$ containing $\fix(\si)(K)$. 
  \end{rem}


\begin{thebibliography}{CHP}
\bibitem{A} Bijan Afshordel, Generic Automorphisms with Prescribed Fixed
  Fields, J. of Symb. Logic 79 (2014), no. 4, 985 -- 1000. 

 \bibitem{ACFAmod} Z. Chatzidakis, A note on canonical bases and one-based types in supersimple
theories, Confluentes Mathematici Vol. 4, No. 3 (2012) 1250004 (34
pages). DOI: 10.1142/S1793744212500041.  
\bibitem{C-Psf} Z. Chatzidakis, A note on the non-existence of prime
  models of theories of pseudo-finite fields. Preprint arXiv 2004.05593.

 \bibitem{CHS} Z. Chatzidakis, C. Hardouin and M. Singer, On the definition of
difference Galois groups,   in: {\em Model Theory with
applications to algebra and analysis, I}  (Z. Chatzidakis,
H.D. Macpherson, A. Pillay, A.J. Wilkie editors), Cambridge University
Press, Cambridge 2008, 73 -- 110.



  
\bibitem{CH} Z. Chatzidakis, E. Hrushovski, Model theory of difference 
fields,
Trans. Amer. Math. Soc. 351 (1999),  2997 -- 3071. 



 




\bibitem{Co} R.M. Cohn, {\em Difference algebra}, Tracts in
Mathematics 17, Interscience Pub. 1965.

\bibitem{Du} J. -L. Duret, Les corps faiblement alg\'ebriquement clos
non s\'eparablement clos ont la propri\'et\'e d'ind\'ependance, in: 
Model theory of Algebra and Arithmetic, Pacholski et al. ed.,
Springer Lecture Notes 834 (1980), 135 --157.

\bibitem{FJ} M. Fried, M. Jarden, {\em Field
Arithmetic},  3rd edition, Ergebnisse 11, Springer Berlin-Heidelberg 2008.

\bibitem{Hr-PAC} E. Hrushovski, Pseudo-finite fields and related
structures, in: {\em Model Theory and Applications}, B\'elair et al. ed.,
Quaderni di Matematica Vol. 11, Aracne, Rome 2005, 151 -- 212.


\bibitem{PS} Anand Pillay and Saharon Shelah, Classification theory over
  a predicate I, Notre Dame J. of Formal Logic
Volume 26, Number 4, October 1985, 361 -- 376.



  \bibitem{La-Ell} Serge Lang, {\em Elliptic Functions}, Graduate Texts
    in Mathematics 112, 2nd edition, Springer-Verlag 1987.

 

    \bibitem{SU} Saharon Shelah and Alexander Usvyatsov, Classification
      over a predicate -- the general case: Part I - structure
      theory. Preprint arXiv:1919.1081.

       \bibitem{Sil} Joseph H. Silverman, {\em The Arithmetic of Elliptic
    Curves}, Grad. Texts in Math. 106, 2nd edition, Springer
    Verlag, 2016. 


\end{thebibliography}
\end{document}